\documentclass[11pt,twoside,a4paper]{article}
\usepackage{amsmath,amssymb,amsthm}

\usepackage{graphicx,psfrag}

\newtheorem{Thm}{Theorem}[section]

\theoremstyle{definition}

\theoremstyle{remark}

\newcommand{\R}{\mathbb{R}}
\newcommand{\Z}{\mathbb{Z}}

\newcommand{\Hy}{\mathbb{H}}

\newcommand{\al}{\alpha}

\newcommand{\ga}{\gamma}

\newcommand{\ep}{\varepsilon}

\renewcommand{\phi}{\varphi}

\newcommand{\diam}{\operatorname{diam}}

\newcommand{\Ric}{\operatorname{Ric}}
\newcommand{\tr}{\operatorname{tr}}

\newcommand{\vol}{\operatorname{vol}}

\newcommand{\eps}{\varepsilon}

\renewcommand{\d}{\partial}

\begin{document}

\title {Almost maximal volume entropy}     
\author{Viktor Schroeder\footnote{Supported by Swiss National
Science Foundation} \ \& Hemangi Shah\footnote{The author thanks 
the Institute of Mathematics of the University of Z{\"u}rich for its hospitality
and support}}
\maketitle

\begin{abstract}
We prove the existence of manifolds with almost maximal volume entropy which are not hyperbolic.
\noindent

\end{abstract}
\section{Introduction}

To a closed Riemannian manifold one associates an asymptotic invariant 
$h_v(M)$ called the {\em volume entropy} by 
$$ h_v(M) = \lim_{r\to \infty} \frac{1}{r} \log \vol(B_r(\tilde{x})), $$
where $B_r(\tilde{x})$ is the ball of radius $r$ about a point $\tilde{x}$ in the 
universal covering space $\tilde{M}$ of $M$. It was shown by Manning \cite{M} that the above
limit always exists and is independent of the choice of $\tilde{x}$.

In \cite{LW} Ledrappier and Wang proved the following remarkable rigidity result
in the case of a lower bound of the Ricci curvature $\Ric$ for $M$:

\begin{Thm}
 Let $M$ be a closed $n$-dimensional Riemannian manifold with
 $\Ric \geq -(n-1)$. Then the volume entropy satisfies $h_v(M) \leq (n-1)$
 with equality if and only if $M$ is hyperbolic, i.e. the curvature is
 constant equal to $-1$.
 
\end{Thm}

In \cite{CRX} Chen, Rong and Xu 
proved the following stability result:

\begin{Thm}
 Given $n$ and $d>0$ there exists $\ep=\ep(n,d)>0$ such that any closed
 $n$-dimensional Riemannian manifold with
 $\Ric \geq -(n-1)$, $\diam(M) \leq d$ and $h_v(M) \geq (n-1) -\ep$ is
 diffeomorphic to a hyperbolic manifold, i.e. it carries a metric of constant
 curvature $-1$.
\end{Thm}

It is natural to ask if this results generalizes to the case, where 
the diameter bound $\diam(M) \leq d $ is replaced by some volume bound $\vol(M) \leq v$.
We give examples that this is not the case. In the following statement $K$ denotes the sectional curvature.

\begin{Thm} \label{Thm:main}
 For every $n \geq 3$ there exists $v > 0$ and a sequence of
 $n$-dimensional closed Riemannian manifolds
 $M_i$ which satisfy the following conditions:
 
 \begin{enumerate}
  \item[(i)] $-1 \leq K \leq 0$. In particular $\Ric \geq -(n-1)$. 
  \item[(ii)] $\vol(M_i) \leq v$.
  \item[(iii)] $h_v(M_i) \geq (n-1) - \ep_i$ with $\ep_i \to 0$.
  \item[(iv)] $M_i$ cannot carry a metric with negative sectional curvature. In particular $M_i$ is not diffeomorphic to a hyperbolic manifold.
 \end{enumerate}

 \end{Thm}
 
 Our examples will be consequences of examples of the following type:
 
 \begin{Thm} \label{Thm:tech}
 For every $n\geq 3$ there exists $v >0$ such that for all $\eps >0$ there exists a closed $n$ dimensional Riemannian manifold $M = M_{\eps}$ with 
 the properties
 \begin{enumerate}
  \item[(i)]  $-1-\eps \leq K \leq 0$. 
  \item[(ii)] $\vol(M)\leq v$.
  \item[(iii)] There exists an open subset $W\subset M$ with $\vol(W) \geq (1-\eps)\vol(M)$, such that $K \equiv -1$ on $W$.
  \item[(iv)] The fundamental group $\pi_1(M)$ has an abelian subgroup of $\operatorname{rank} \geq 2$.
 \end{enumerate}

 \end{Thm}

 We will give two types of constructions which lead in a similar way to the examples of Theorem \ref{Thm:tech}. One is 
 a doubling construction (which works in all dimensions $\geq 3$), the other is a cusp closing construction which 
 give examples in dimension $\geq 4$.  
 
 We show in section \ref{sec:pomt} how the main Theorem \ref{Thm:main} follows from \ref{Thm:tech}.
 In section \ref{sec:ccl} we give the details of the cusp closing construction and in section \ref{sec:doub} we discuss the doubling construction.
 In the final section we pose a question regarding the possibility of cusp closings under certain curvature conditions.
 
 The main part of the work was done during a visit of the second author at the University of Zurich in November and December 2016. 
 The first author thanks the University of Oxford and the Newton Institute in Cambridge for the hospitality.
 It is also a pleasure to thank 
 Manfred Einsiedler, Ruth Kellerhals, Camillo de Lellis and Alain Reid for inspiring discussions and email exchanges.
 
 \section{Entropies} \label{sec:entr}
 
 In this paper we use three versions of entropy. First the {\em volume entropy} $h_v$ which was already explained in the 
 introduction. 
 
 The {\em topological entropy} $h_t$ can be defined for a pair $(X,T)$, where $X$ is a compact
 topological space and $T$ is a homeomorphism of $X$. The original definition is due to Adler, Konheim and McAndrew \cite{AKM}.
 Three different but equivalent definitions are given in \cite{M}. 
 In the case of a closed Riemannian manifold $M$ we understand under $h_t(M)$ the topological entropy of the pair
 $(SM,\Phi_1)$, where $SM$ is the unit tangent bundle and $\Phi_1$ is the time $1$ map of the geodesic flow.
 In his seminal paper \cite{M} Manning proved that for all closed Riemannian manifolds
 $h_t(M)\geq h_v(M)$. He further showed that $h_t(M)=h_v(M)$ in the case that all sectional curvatures are nonpositive.
 This equality case was generalised by Freire and Ma\~n\'e \cite{FM} to the case that $M$ has no conjugate points.
 
 The {\em measure theoretic entropy} $h_{\mu}$ was introduced by Kolmogorov \cite{K} and Sinai \cite{Si}
 and is defined for a measure preserving map $T$ on a probability space $(X,\cal B,\mu)$.
 In the case of a closed Riemannian manifold we view $h_{\mu}(M)$ as the 
 measure theoretic entropy of the time $1$ map of the geodesic flow on $(SM,\cal B,\mu)$, where $\cal B$ is the Borel $\Sigma$-algebra
 on the unit tangent bundle and $\mu$ is the normalized Liouville measure on $SM$.
 
 It was shown by Goodwyn \cite{G1} that in the case that $T$ is a homeomorphism of a metric space 
 $(X,d)$, the topological entropy $h_t$ is greater or equal $h_{\mu}$ for any $T$-invariant probability measure 
 $\mu$ on the Borel sets of $(X,d)$. Later he showed \cite{G2} that $h_t = \sup_{\mu} h_{\mu}$, where $\mu$ is 
 a $T$-invariant probability measure.
 In particular $h_t(M) \geq h_{\mu}(M)$ for a closed Riemannian manifold $M$.

 \section{Proof of the main Theorem} \label{sec:pomt}
 
 We show in this section that Theorem \ref{Thm:tech} implies our main result Theorem \ref{Thm:main}.
 
 Note that by the theorem of Preissmann \cite{P} the resulting manifolds $M_{\eps}$ cannot carry a metric of negative curvature,
 since
 the fundamental group contains an abelian subgroup of rank $\geq 2$. 
 
 We use a result of Ballmann and Wojtkowski \cite{BW} to estimate the measure theoretic entropy $h_{\mu}(M)$ of the geodesic flow
 on $SM$, where $\mu$ is the normalized Liouville measure on the unit tangent bundle $SM$. 
 For a unit tangent vector
 $v \in S_pM$ consider the Jacobi operator $R(\cdot,v)v$, which is a nonpositive symmetric operator on $T_pM$.
 The formula in \cite{BW} gives
 $$h_{\mu}(M) \geq \int_{SM} \tr \sqrt{-R(\cdot,v)v} d\mu(v). $$
 Since $\tr \sqrt{-R(\cdot,v)v)} = (n-1)$ for $v\in T_pM$, $p\in W$ we obtain
 $h_{\mu}(M_{\eps}) \geq (n-1)(1-\eps)$.
 If we rescale the metric on $M_{\eps}$ such that the new lower curvature bound is  $-1$ we obtain
 closed Riemannian manifolds $\bar{M}_{\eps}$ with $-1 \leq K \leq 0$ and $h_{\mu} (\bar{M}_{\eps}) \geq (n-1) -\bar{ \eps}$ with $\bar{\eps}\to 0$ if $\eps\to 0$.
 Now Goodwyn's result ($h_t \geq h_{\mu}$) and Manning's result ($h_t=h_v$, since $K\leq 0$) shows
 $h_{\mu} \leq h_v$. Thus we obtain the desired examples.

 \section{The cusp closing construction} \label{sec:ccl}

 We first describe informally the idea of the cusp closing construction. We start with a complete open $n$-dimensional hyperbolic
 (constant curvature $-1$) Riemannian manifold of finite volume. Such a manifold has finitely many 
 ends which are cusps. Thus they are isometric to
 $\operatorname{HB}/\Gamma$, where $\operatorname{HB}$ is a horoball in $\Hy^n$ around some point
 $z \in \partial \Hy^n$ and $\Gamma$ is a group of parabolic isometries fixing $z$ and operating as a lattice on
 the horosphere $\operatorname{HS}=\partial \operatorname{HB} \simeq \R^{n-1}$.
 
 We assume that all the cusp cross section are tori, i.e. we have $\Gamma \simeq \Z^{n-1}$.
 Hyperbolic manifolds of such type exists in every dimension by a result of McReynolds, Reid and Stover \cite{MRS}.
 
 Let now $\eps > 0$ be given. Then we will cut off all cusps very far out such that the cut off parts have only
 $\eps$ amount of volume. After cutting we obtain a compact manifold with finitely many boundary components,
 which are all flat tori in their intrinsic horospherical metric and all of this tori have very small volume.
 
 Now we close the cusp with a tube around a closed torus of codimension $n-2$. Such a construction was
 carried out in \cite{S}.
 We also have to control the curvature and the volume of the cusp closing.
 Unfortunately we are not able to close the cusp with curvature bound $\Ric \geq -(n-1)$ and $K \leq 0$.
 Actually we conjecture that a cusp closing is impossible with this curvature restrictions (see section \ref{sec:quest}).
 However, we can close the cusp such that the curvature satisfies $-1-\eps \leq K \leq 0$ and the volume
 of each part is also of size $\eps$. 
 
 Using this construction we obtain a compact Riemannian manifold with two parts:
 A {\em thick} part which contains almost all of the volume and where the curvature is constant $-1$.
 The remainder of the manifold consists of 
 a finite number of tubes around closed tori. These tubes have very large diameter but tiny volume and
 sectional curvature bounded by $\-1-\eps \leq K \leq 0$.

 We now give the details of the construction.
 As above a cusp is written as $\operatorname{HB}/\Gamma$, and
 we can write $\operatorname{HB}$ metrically as warped product
 $(0,\infty)\times_{f} \R^{n-1}$, with the metric
 $dt^2 + f^2(t) d\sigma^2$, where $d\sigma$ is the flat metric on
 $\R^{n-1}$ and $f(t)= e^{-t}$. The goup $\Gamma \simeq \Z^{n-1}$ operates
 as a lattice on $\R^{n-1}$.
 
 We want to fill the cut off cusp with a tube around a torus of
 dimension $n-2$. To describe the universal cover of the tube we consider cylindrical coordinates on
 $\R^{n}$, i.e. we consider
 $[0,\infty) \times S^1 \times R^{n-2}$ with a metric of type
 $dt^2 + s^2(t) d\varphi^2 + c^2(t) d \sigma^2$, where
 $dt^2, d\varphi^2$ and $d\sigma^2$ are the standard metrics on the factors.
 In order that the metric is smooth at $t=0$ we consider smooth functions
 $s,c:[0,\infty)\to [0,\infty)$ with $s(0)=0, s'(0)=1$, $c(0)=1, c'(0)= 0$, such that
 $s$ extends as an odd and $c$ as an even function to $\R$.
 Note that in the case $n=3$ and functions $c=\cosh$ and $s=\sinh$ we obtain a metric description
 of hyperbolic $3$-space in cylindrical coordinates.
 
 By the standard warped product formulas we can compute the curvature
 tensor and the relevant sectional curvatures.
 If we consider the tangent vectors $\frac{\d}{\d t}$, $\frac{\d}{\d \varphi}$ and tangent vectors
 $U,V$ to $\R^{n-2}$ ( $U$, $V$ linearly independent), then the sectional curvature of the corresponding 
 planes are
 
   $$ K(\frac{\d}{\d t},\frac{\d}{\d \varphi})=-\frac{s''}{s}, \ \ \ \ \ \
  K(\frac{\d}{\d t},U)=-\frac{c''}{c}, $$
  $$ K(\frac{\d}{\d \varphi},U)=-\frac{s'c'}{sc},\ \ \ \ \ \
   K(U,V)= -\frac{c'^2}{c^2}. $$

If we are in dimension $\geq 4$ and we choose $c=\cosh$ and $s=\sinh$ we obtain a metric  
on $\R^n$, which contains a totally geodesic and flat $\R^{n-2}$ (corresonding to $t=0$) and which has negative
curvature away from this flat and if $t \to \infty$ the metric tends to a metric of constant curvature $\equiv -1$.
In order that the metrics finally match up with the cut off cusp, we will choose the functions
 $s$ and $c$ in a way that for $t$ sufficiently large (say $t\geq r_0$) we have
 $s'(t) = s(t)$ and $c'(t)=c(t)$, i.e. $s(t)=c_1 e^t$ and $c(t)= c_2 e^t$.
 Actually we will construct the functions in a way that
 $s(t)=c(t)=\frac{1}{2}e^t$ for $t\geq r_0$.
 This implies that for $t_0 \geq r_0 + 1$ the collar
 $(t_0 -1 ,t_0)\times_{s}S^1\times_{c}\R^{n-2}$ has constant curvature $-1$ and
 is isometric to $((0,1)\times_{f}\R^{n-1})/<\gamma>$ where $<\gamma>$ is a cyclic group of isometries
 generated by a translation $\gamma$ on $\R^{n-1}$ with translation length
 $2\pi s(t_0)$. Note that under this identification
 $\{t_0-1\}\times S^1 \times \R^{n-2} $ corresponds to $\{1\}\times \R^{n-1} / <\gamma>$ and
 $\{t_0\}\times S^1\times \R^{n-2} $  to $\{0\}\times \R^{n-1} / <\gamma>$ .
  
More concretely for given
$\eps >0$ we define functions $c_{\eps}$ and $s_{\eps}$ by
$$2\ c_{\eps}(t)= e^t+\varphi_{\eps}(t) e^{-t}$$
$$2\ s_{\eps}(t)= e^t-\varphi_{\eps}(t) e^{-t}$$
where $\varphi_{\eps}$ is a suitable cut off function with
$0\leq \varphi_{\eps} \leq 1$, which is constant equal $1$ on $[0,1]$,
monotonously nonincreasing, constant equal to $0$ for $t\geq r_{\eps}$, such that first and second
derivatives are small. 
Thus for $t\leq 1$ the metric coincides with the case $c=\cosh$, $s=\sinh$ and for $t\geq r_{\eps}$
the metric has constant curvature $\equiv -1$.
Using the above warped product curvature formulas one sees that we can choose the function
$\varphi_{\eps}$ in a way that 
for the sectional curvature $-1-\eps \leq K \leq 0$. Of course $r_{\eps} \to \infty$ if $\eps \to 0$.

We now have to arrange that the corresponding group actions fit together.
The cusp has a collar isometric to
$((0,1)\times_{f}\R^{n-1})/\Gamma$, where $\Gamma$ is a discrete cocompact group of translations 
of $\R^{n-1}$, which we can identify with a subset of $\R^{n-1}$.
We choose a set $\al_1,\ldots,\al_{n-1}$ of generators of $\Gamma$ in the following standard way:
let $\al_1 \in \Gamma \setminus \{0\}$ of minimal length, then
$\al_2 \in \Gamma\setminus <\al_1>$ of minimal length etc.
In suitable coordinates we have
$\al_1 =(a_{11},0,\ldots,0) \in \R^{n-1}$, $\al_2=(a_{21},a_{22},0,\ldots,0)$ etc.

For any integer $k \geq 1$ also 
$\ga_1:= k \al_1 + \al_2$ , $\ga_2:= \al_2$ ...$ \ga_{n-1}:= \al_{n-1}$ is a system of generators.
We choose now $k$ large enough such that
the translation length
$\|\ga_1\|=\|k\al_1 + \al_2\| \geq 2\pi s_{\eps}(r_{\eps}+1)$.
Thus we can find $t_0 \geq r_{\eps}+1$ such that $\|\ga_1\|=2\pi s_{\eps}(t_0)$.
This implies by the above arguments that 
$((0,1) \times_{f} \R^{n-1} )/<\ga_1>$ is isometric to 
$(t_0-1,t_0)\times_{s_{\eps}}S^1\times_{c_{\eps}}\R^{n-2}$.
The whole collar $((0,1)\times_f\R^{n-1})/\Gamma$ can be written as
$((0,1)\times_f\R^{n-1})/<\ga_1>\times \Delta$, where
$\Delta =<\al_2,\ldots,\al_{n-1}> \simeq \Z^{n-2}$.
Note that the operation of $\Delta$ can be isometrically extended to an operation on
$[0,t_0+1)\times_{s_{\eps}}S^1\times_{c_{\eps}}\R^{n-2}$ and hence we can isometrically glue
$([0,t_0+1)\times_{s_{\eps}}S^1\times_{c_{\eps}}\R^{n-2}) / \Delta$ to the cusp.

Finally one easily estimates that the $n$ dimensional volume of the tube 
$([0,t_0+1)\times_{s_{\eps}} S^1\times_{c_{\eps}} \R^{n-2}) / \Delta$
is estimated by the $n-1$ dimensional volume of the boundary torus $\R^{n-1}/ \Gamma$, since $s_{\eps}$ is exponential.
Since the boundary tori have small volume also the tubes have small volume.

\section{The doubling construction} \label{sec:doub}

We describe the idea of the doubling construction.
The idea of this construction goes back to Heintze. Again we start with a complete open
$n$-dimensional hyperbolic manifold $M$ of finite volume which has (assume for simplicity) one cusp
isometric to $((0,\infty)\times_f\R^{n-1})/\Gamma $.

Now we change the warping function $f$ to a function which is still strictly convex but which has a unique
positive minimum. Thus consider a symmetric function $c:\R \to (0,\infty)$ , $c(t)=c(-t)$, $c''(t) > 0$ such that 
$c'(t)=c(t)$ for $t \geq t_0$.
Then $(-t_0-1,t_0+1)\times_{c} \R^{n-1}$ has two mirror symmetric collars
$(-t_0-1,-t_0)\times_c \R^{n-1}$ and $(t_0,t_0+1)\times_c \R^{n-1}$which are both
isometric to $(0,1)\times_f \R^{n-1}$ and we can arrange the operation of the group $\Gamma$,
such that 
$(-t_0-1,t_0+1)\times_{c} \R^{n-1} / \Gamma$ can be glued isometrically on both sides onto two copies 
of the manifold minus the cusp. For a suitable $c = c_{\eps}$ one can show that the resulting manifold has curvature
$-1-\eps \leq K \leq 0$ and that also the volume of the channel between the two copies of the original manifold is small.
The torus $\{0\}\times \R^{n-1}/\Gamma$ is a flat totally geodesic torus.
In the case of more cusps we glue two copies of the manifold simultaneously along each pair of cusps.

\section{Question}  \label{sec:quest}

To our opinion the construction has an esthetical flaw, namely the curvature
condition $-1-\eps \leq K $ in point (i) of Theorem \ref{Thm:tech} . The $\eps$ in this point
makes the {\em ugly} rescaling in section \ref{sec:pomt} necessary.

We tried to obtain a construction which instead of the curvature condition $-1-\eps \leq K \leq 0$ satisfies

\begin{equation} \label{eq:1}
 -(n-1)\leq \Ric, \ \ K\leq 0.
\end{equation}

However we have arguments that the doubling construction does not work with the condition (\ref{eq:1}).
In the cusp closure construction one has much more flexibility due to the two functions $s$ and $c$.
But finally we think that a cusp closing under curvature condition
$-(n-1)\leq \Ric, K\leq 0$ is not possible, which would be a surprising rigidity result.
Actually we conjecture that the following question has a positive answer:

Let $M$ be a closed $n$-dimensional Riemannian manifold with $-(n-1) \leq \Ric$ and $K\leq 0$. Let
$M = M_{\geq \mu} \cup M_{\leq \mu}$ be the thick-thin decomposition, where $M_{\geq \mu}$ (resp. $M_{\leq \mu}$) is the part
where the injectivity radius is larger (resp. smaller) than the corresponding Margulis constant $\mu$.
Assume $K \equiv -1$ on $M_{\geq \mu}$.

Is then every abelian subgroup of $\pi_1(M)$ cyclic ( i.e. is there no
higher rank abelian subgroup) ?

Institut f{\"u}r Mathematik, Universit{\"a}t Z{\"u}rich, 
Winterthurerstrasse 190,\\
\indent
CH-8057 Z{\"u}rich.\\
\indent
Email :  viktor.schroeder@math.uzh.ch\\\\
\indent
Harish Chandra Research Institute, HBNI,  Chhatnag Road,\\
\indent
Jhusi, Allahabad 211019, India.\\
\indent
Email :  hemangimshah@hri.res.in

\end{document}